\documentclass[aap,MSNbibl,dvips]{arximspdf}

\doi{10.1214/12-AAP911} %kopijuoti is PTS
\volume{24}
\issue{2}
\pubyear{2014}
\firstpage{449}
\lastpage{475}

\makeatletter

\newcommand{\rrvert}{\vert}
\newcommand{\llvert}{\vert}

\newcommand{\indic}[1]{\mathbf{1}_{#1}}
\newcommand{\indica}[1]{\mathbf{1}_{\{#1\}}}

\newtheorem{theorem}{Theorem}

\newtheorem{lemma}[theorem]{Lemma}
\newtheorem{coro}[theorem]{Corollary}
\newtheorem{proposition}[theorem]{Proposition}
\newproclaim{defn}[theorem]{Definition}
\newproclaim{remark}[theorem]{Remark}

\newcommand{\eqref}[1]{(\ref{#1})}

\newcommand{\eps}{\varepsilon}
\newcommand{\E}{\mathbb{E}}
\renewcommand{\P}{\mathbb{P}}
\newcommand{\R}{\mathbb{R}}
\newcommand{\Q}{\mathbb{Q}}
\newcommand{\bck}{}
\newcommand{\La}{\Lambda}
\newcommand{\sfrac}[2]{{\frac{#1}{#2}}}
\newcommand{\N} {\mathbb{N}}
\newcommand{\eqlaw}{\stackrel{d}{=}}

\def\cR{\mathcal{R}}

\makeatother

\begin{document}
\begin{frontmatter}

\title{A small-time coupling between $\Lambda$-coalescents and
branching processes}
\runtitle{A small-time coupling between $\Lambda$-coalescents and CSBPs}

\begin{aug}
\author[A]{\fnms{Julien} \snm{Berestycki}\corref{}\thanksref{t1}\ead[label=e1]{julien.berestycki@upmc.fr}},
\author[B]{\fnms{Nathana\"{e}l} \snm{Berestycki}\thanksref{t2}}
\and\break
\author[C]{\fnms{Vlada} \snm{Limic}\thanksref{t3}}
\thankstext{t1}{Supported in part by ANR A3 and ANR MADCOF.}
\thankstext{t2}{Supported in part by EPSRC Grants EP/GO55068/1 and EP/I03372X/1.}
\thankstext{t3}{Supported in part by NSERC Discovery grant, by Alfred P. Sloan Research Fellowship and by ANR MAEV research grant.}
\runauthor{J. Berestycki, N. Berestycki and V. Limic}
\affiliation{Universit\'e Pierre et Marie Curie---Paris VI, University
of Cambridge and Universit\'{e} Paris Sud 11}
\address[A]{J. Berestycki\\
LPMA, UMR 7599\\
Universit\'e Pierre et Marie Curie---Paris VI\\
Sorbonne Universit\'{e}s\\
F-75005, Paris\\
France\\
\printead{e1}} %adresu isvedimo komanda gale!
\address[B]{N. Berestycki\\
CMS\\
University of Cambridge\\
Wilberforce Rd.\\
Cambridge CB3 0WB\\
United Kingdom}
\address[C]{V. Limic\\
UMR 8628\\
D\'{e}partement de Math\'{e}matiques\\
Universit\'{e} Paris Sud 11\\
91405 Orsay\\
France}
\end{aug}

% HISTORY:
\received{\smonth{6} \syear{2012}}
\revised{\smonth{11} \syear{2012}}

% ABSTRACT
%
\begin{abstract}
We describe a new general connection between
$\Lambda$-coalescents and genealogies of continuous-state branching processes.
This connection is based on the construction of an explicit coupling
using a particle representation inspired by the lookdown process of
Donnelly and Kurtz. This coupling has the property that the coalescent
comes down from infinity if and only if the branching process becomes
extinct, thereby answering a question of Bertoin and Le Gall.
The coupling also offers new perspective on the speed of coming down
from infinity and allows us to relate power-law behavior for $N^\La
(t)$ to the classical upper and lower indices arising in the study of
pathwise properties of L\'evy processes.
%This allows us to give an intuitive explanation for the sampling
%formulae and to revisit a number of recent results.
\end{abstract}

% KEYWORDS
% Pirmas kwd is didziosios raides
%
\pdfkeywords{60J25, 60F99, 92D25, $\Lambda$-coalescents, continuous-state branching processes, particle system representation,
Levy processes, lookdown construction, Fleming--Viot processes}

\end{frontmatter}

\section{Introduction and main results}\label{sec1}
Coalescents with multiple collisions,
also known as $\Lambda$-\emph{coalescents} are Markovian models of
coagulation. Introduced and first studied independently by Pitman \cite
{pit99} and by Sagitov \cite{sag99} (also considered in a
contemporaneous work of Donnelly and Kurtz \cite{dk99}), these
processes have been intensely researched in the last decade.
The research is mostly motivated by the fact that $\La$-coalescents
arise naturally as scaling limits for the genealogy in certain
exchangeable population dynamics models.
We refer to \cite{bertoin,ensaios} for an introduction and a survey
of the relevant literature.

The \textit{standard} $\Lambda$-coalescent starts with infinitely many
microscopic particles that coalesce into larger clusters as time runs.
Our interest in this paper concerns the small-time behavior of
(standard) $\La$-coalescents, in particular the phenomenon of \emph
{coming down from infinity} (a precise definition will be given below).
Our main goal is to answer a question which arose from work of Bertoin
and Le Gall \cite{blg3}. They observed that the Schweinsberg condition
\cite{sch1} for coming down from infinity for $\Lambda$-coalescents
is equivalent to the condition for extinction of related
continuous-state branching processes (CSBPs),\vadjust{\goodbreak} and asked if a deeper
connection exists between these two classes of processes.
In this paper, we construct an explicit coupling between a given
$\Lambda$-coalescent and a certain associated CSBP, and therefore
answer the above question of Bertoin and Le Gall.

This coupling makes use of a particle system representation based on a
lookdown process in the spirit of Donnelly and Kurtz \cite{dkfirst,dk99}. %This approach turns out to be the best way to
Apart from its interest from a purely theoretical point of view, our
coupling gives a new understanding of the asymptotic form of the
``speed of coming down from infinity'' (as discussed by the authors in
\cite{bbl1}) and leads to precise quantitative results for the
corresponding $\La$-coalescent observed at small times. In particular,
the power-law exponents for the number of blocks in a particular $\La
$-coalescent are shown to coincide with the classical notion of upper
and lower indices of the L\'evy measure of the associated CSBP.

%The reader might also be interested to know the limitations of this
%technique as explained in Remark~\ref{Rlimtech}.
The methodology in this paper has several points in common with \cite
{bbs1,bbs2}, where an analogous link between Beta-coalescents and
$\alpha$-stable continuous-state branching processes was used.
However, in these papers the central tool was an explicit embedding of
the lookdown process into the (stable) continuous random tree, which
allowed for many explicit computations. Here, we show that the correct
way to generalize this picture for an arbitrary $\Lambda$-coalescents
is directly via the particle system approach of the lookdown process.

In the rest of the paper, we denote by $\eqlaw$ the equivalence in
distribution.
We also use the standard Bachmann--Landau notation $\sim, O(\cdot),
o(\cdot),\asymp$ for comparing asymptotic behavior of deterministic
and stochastic functions and sequences.

\subsection{Coalescents and CSBPs}\label{sec1.1}

Let $\Lambda$ be an arbitrary finite measure on $[0,1]$, and let $(\Pi
_t, t \ge0)$ denote the associated $\Lambda$-coalescent. The Markov
jump process $(\Pi_t, t \ge0)$ takes values in the set of partitions
of $\{1, 2, \ldots\}$. Its law is specified by the requirement that,
for any $n\in\N$, the restriction $\Pi^n$ of $\Pi$ to $\{1, \ldots, n\}$ is a continuous-time Markov chain with transition rates given as
follows: whenever $\Pi^n$ has $ b \in[2, n]$ blocks, any given
$k$-tuple of blocks coalesces at rate $\lambda_{b,k}:= \int_{(0,1]}
r^{k-2}(1-r)^{b-k} \La(dr)$.\vspace*{1pt}

We will always assume that $\Pi(0)$ is the trivial partition $\{\{i\}
\dvtx i\in\N\}$. Let us call $N^{\La}(t)$ the number of blocks of $\Pi
(t)$ the coalescent at time $t$. The first question one may ask about
these processes is whether the number of blocks ever becomes finite. In
his seminal paper \cite{pit99} Pitman noted that [provided $\La(\{1\}
)=0$] as a consequence of the strong Markov property, the following
striking dichotomy holds: either $\P(N^{\La}(t) =\infty, \forall
t\ge0)=1$ or $\P(N^{\La}(t) <\infty, \forall t > 0)=1$. In the
latter case the coalescent is said to \emph{come down from infinity}.
Finding a necessary and sufficient condition for this phenomenon was
naturally one of the first problems to be studied. As part of his
thesis work, Schweinsberg \cite{sch1} derived the following criterion:
the $\La$-coalescent comes down from infinity if and only if
\begin{equation}
\label{Econdsch} \sum_{b=2}^\infty \Biggl( \sum
_{k=2}^b (k-1) \pmatrix{b \cr k} \lambda
_{b,k} \Biggr)^{-1} <\infty.
\end{equation}

Over the subsequent years, a series of remarkable links were discovered
between $\Lambda$-coalescents and continuous-state branching processes
(CSBP), for some special cases of $\Lambda$. The case of Kingman's
coalescent ($\Lambda= \delta_{0}$) was analyzed by Perkins \cite
{perkins} in 1991, though he used a somewhat different language.
Bertoin and Le Gall \cite{blg0} studied the case of the
Bolthausen--Sznitman coalescent (where $\Lambda(dx) = dx$ is the
uniform measure on $[0,1]$), and then Birkner et al. \cite{7} studied
all the Beta-coalescents cases
[where $\La$ is the $\operatorname{Beta}(2-\alpha, \alpha)$ distribution, and
$\alpha\in(0,2)$].

While seeking a way to understand the above results as special cases of
a general theorem, Bertoin and Le Gall \cite{blg3} made the following
observation. Consider the function
\begin{equation}
\label{Dpsi} \psi(q):= \int_0^1
\bigl(e^{-qx} -1+qx \bigr)x^{-2}\La(dx), \qquad  q \ge0.
\end{equation}
Then $\psi$ is the Laplace exponent of a spectrally positive L\'evy
process and is thus the \emph{branching mechanism} of a CSBP $(Z_t,
t\ge0)$. (Definitions and elementary properties of CSBPs may be found,
for instance, in \cite{lamp,lamperti2} and \cite{ensaios} and
Chapter~6 of \cite{athreya-ney}.) In particular Grey \cite{grey}
showed that a $\psi$-CSBP becomes extinct almost surely in finite time
if and only if
\begin{equation}
\label{Egrey0} \int_1^\infty{dq \over\psi(q) }
<\infty.
\end{equation}

\begin{theorem}[(Bertoin and Le Gall, \cite{blg3})] \label{Tmainthm}
Conditions \eqref{Econdsch} and \eqref{Egrey0} are equivalent. In
other words, a particular (standard) $\La$-coalescent comes down from
infinity if and only if the corresponding CSBP becomes extinct.
\end{theorem}
The proof of Bertoin and Le Gall (see the end of Section~4 in \cite
{blg3}) is direct and analytical. However, Theorem~\ref{Tmainthm}
strongly suggests that a general probabilistic connection exists, and
this prompted Bertoin and Le Gall to ask for a probabilistic proof of
their result.

%Let an arbitrary finite measure $\Lambda$ be given.
The main goal of the present work is to provide an explicit coupling
that makes Theorem~\ref{Tmainthm} ``obvious.''\vadjust{\goodbreak} In fact, the coupling
yields much more information, including a quantitative estimate on
$N^\La(t)$ for small times $t$ (Propositions \ref{Tsmalltime} and
\ref{Palmostsure}). This estimate matches the ``speed of coming down
from infinity'' obtained by the authors in \cite{bbl1} with a
martingale method. In fact, the present coupling construction suggested
that completely general result in the first place.

\textit{Organization and contents of the paper}.
Our coupling is based on a particle system representation for $\La
$-coalescents and a connection to a version of Donnelly and Kurtz's
\emph{lookdown process.} Both for the sake of completeness and of
explaining the
differences between our construction and that of \cite{dk99}, we will
start by defining the lookdown process. More precisely, we will show
that this construction is feasible whenever its driving point process
$\pi= \sum_i \delta_{(t_i, p_i)}$, given on $(0, \infty) \times
(0,1)$ satisfies $\sum_{t_i \le t } p_i^2 < \infty$ for all $t\ge0$.
This result, which we believe is of independent interest, is stated in
Proposition~\ref{Plookdowngeneral}.

We will then apply this construction to two distinct point processes,
one arising from the $\La$-coalescent and the other from the
associated CSBP. This is done in Section~\ref{Spartsystems}. We then
use these representations to obtain a coupling between the two
processes. This allows us to conclude that the genealogy of the CSBP
is, at small times, ``close'' to the $\La$-coalescent. On the other
hand, the CSBP gets extinct in finite time if and only if the number of
individuals with descendants alive at a future time $t>0$ is finite
(Proposition~\ref{PNXandNZ}).
This directly yields Theorem~\ref{Tmainthm} and its stronger
quantitative version, Theorem~\ref{Tsmalltime}.

We next use these results together with certain pathwise properties of
L\'evy processes and CSBPs to discuss the regularity of $N^\La(t)$ as
$t \to0$. Our main result there (Proposition~\ref{partial2}) shows
that the power-law behavior for $N^\La(t)$ is intimately related to
the classical \emph{upper and lower indices} of the L\'evy measure of
$\psi$, following Blumenthal and Getoor \cite{BG} and Pruitt \cite{pruitt}.
The \hyperref[Sexample]{Appendix} contains an example of a measure $\La$ that is not
``well-behaved,''
in the sense that the corresponding $\La$-coalescent comes down from infinity
but the lower and the upper indices are different. We show how this
leads to truly oscillatory behavior for $N^\La(t)$, which highlights
potential difficulties in the analysis of small-time behavior of
general $\La$-coalescents.

\section{Preliminaries}\label{sec2}
\label{Spartsystems}

In this section we describe a general procedure known as the
\textit{lookdown construction}, enabling one to construct measure valued processes
from point processes on $[0,1] \times\R_+$. The material discussed in
this section is mostly well known, but we prefer to give a
brief account of the theory to set the ground for the construction of
the coupling in Section~\ref{Scoupling}.
%Throughout of this section we assume that $\La(\{0\})=0$.
Unless stated otherwise, we henceforth assume that $\La(\{0\})=\La(\{
1\})=0$.

\subsection{Lookdown construction}\label{sec2.1}

The lookdown construction was first introduced by Donnelly and Kurtz in
1996 \cite{dkfirst}.
Their goal was to give a construction of the Fleming--Viot superprocess
that provides an explicit description of the genealogy
of the individuals in the population; see \cite{etheridge} for a
reader-friendly introduction to these notions. Donnelly
and Kurtz subsequently modified their construction in \cite{dk99} to
include more general measure-valued processes
(such as the Dawson--Watanabe superprocesses). It is this version that
we use here, and that we will apply to the
generalized Fleming--Viot superprocesses (which are dual to $\La
$-coalescents) as well as to the ratio processes
associated to CSBPs. Our approach here shares common points with that
of \cite{7}.

For a given (infinite size) population evolving in continuous time, let
the genetic types of individuals be encoded as numbers in $[0,1]$.
More precisely, for each $i \ge1$ and $t\geq0$, let $\xi_i(t)\in[0,1]$
be the genetic type of the individual $i$ (or level~$i$) at time $t$.
As will be seen soon, for our models,
the infinite particle system $((\xi_1(t),\xi_2(t), \ldots),  t\ge0)$
is such that the limiting empirical measure
\[
\Xi_t(\cdot) = \lim_{n \to\infty} \frac{1}n \sum
_{i=1}^n \delta _{\xi_i(t)}(\cdot)
\]
exists simultaneously for all $t$, almost surely.
The process $(\Xi_t(\cdot), t\geq0)$ is a convenient way to track
the evolution of the genetic composition of the population.

We first offer an informal description followed by a formal one in
Definition~\ref{Dlookdownlabel}.
The evolution of $(\xi_i(t))_{i\ge1}$ is driven by a point process
(i.e., a~countable collection of random points) $\pi=(p_i,t_i)_{i\in
\N}$ in $[0,1] \times\R_+$,
and a family of i.i.d. coin tosses.
Each atom of $\pi$ corresponds to a \textit{birth} (or resampling) event.
Changes in $(\xi_i(t), t\ge0)_{i\ge1}$ occur only at birth event times.
Let $(p,t) \in\pi$. Then at time $t$, for each
level $i\geq1$, a coin is tossed, where the probability of head equals
$p$, independently over levels.
Those levels for which the coin comes up heads
(let us denote this set by $I_{p,t}$)
%=\{i_1 < i_2 <...\}$)
modify their label to $\xi_{\min I_{p,t}} (t-)$. In words, each level
in $I_{p,t}$ immediately adopts
the type of the smallest level participating in this birth event.
For the remaining levels reassign the types so that their relative
order immediately prior to this birth event is preserved.
More precisely, for each $i \notin I_{p,t}$, let $\xi_i(t-) = \xi
_{\phi(i)}(t)$ where $\phi$ is
the unique increasing bijection from $\N\setminus\{\min I_{p,t}\}$
onto $\N\setminus I_{p,t}$.

%The coupling between $\Lambda$-coalescents and CSBP that we construct
%is based on applying this construction with two differently
%distributed point processes $\pi$. It is therefore useful to consider
%a more general framework as follows.
A more formal description follows.
Fix $(U_{i,j})_{i,j \ge1}$, a collection of
i.i.d. uniform variables on $[0,1]$. Let $\pi=\{(p_i,t_i)\dvtx i \in\N\}
$ be a \emph{fixed} point process on $[0,1] \times\R_+$ such that
for any $0\le t < \infty$,
\begin{eqnarray}
\label{EconditionH} \sum_{i: t_i\le t} p_i^2<
\infty. \end{eqnarray}
[When we apply this construction later, $\pi$ will be random and we
will work conditionally given $\pi$. Condition \eqref{EconditionH}
will then hold almost surely.]
For each $n\ge1$, construct the label process associated with $\pi$
as follows.
We fix an infinite sequence of exchangeable random variables $(\xi
_i(0))_{i \ge1}$. Set $\xi_i^n(0)=\xi_i(0)$, $i=1,\ldots,n$.
For each $j\ge1$ and $i \in\{1,\ldots,n\}$ define
\begin{eqnarray}
\label{EAiota} A_i(t_j,p_j)\equiv
A_j(i)&:=&\{U_{i,j}\leq p_j\} \quad\mbox{and}
\nonumber
\\[-8pt]
\\[-8pt]
\nonumber
i_1(j)&:=&\min\bigl\{i \geq1\dvtx A_j(i) \mbox{ occurs}\bigr
\}.
\end{eqnarray}
%
%and for all $j\ge1$,
% i_1(j):=\min\{i \geq1: A_j(i) \mbox{ occurs}\}.
For $i\leq n$, let
\begin{equation}
\label{Emiota} m_j(i):=\sum_{l=1}^{i}
\indic{A_j(l)},\qquad i\geq1,
\end{equation}
be the number of levels smaller or equal to $i$ that participate in the
birth event $(p_j,t_j)$.
Denote by $J$ the set of atom indices $\{j\ge1\dvtx m_j(n) \ge2\}$
for which two or more levels in $\{1,\ldots,n\}$ participate in the
corresponding birth event.
Order the collection of indices in $J$ so that $t_{j_1}< t_{j_2}<
\cdots$; this is almost surely possible due to~(\ref{EconditionH});
see Proposition~\ref{Plookdowngeneral}, below.
Define $(\xi_i^n(t))_{1\le i \le n}$ to be constant over $[t_{j_k},
t_{j_{k+1}})$.
Moreover, if $j \in J$,
modify the labels at time $t_j$ as follows: for each $1\le i \le n$ declare
\begin{equation}
\label{Epushup} \xi_i^n(t_j)=
\xi_{ i- (m_j(i)-1)_+}^n(t_j-)\indic {A_j(i)^c}+
\xi_{i_1(j)}^n(t_j-)\indic{A_j(i)},
\end{equation}
where $m_j(i)$ is defined in (\ref{Emiota}).

Finally, observe a crucial property of the above construction:
if $1\le m<n$, then
the restriction of $\xi^n$ to the first $m$ levels yields
$\xi^m$, and in symbols,
\begin{equation}
\label{Ecomprestri} \bigl(\bigl(\xi_1^n(t),\ldots,
\xi_m^n(t)\bigr), t\geq0\bigr)\equiv \bigl(\bigl(
\xi_1^m(t),\ldots,\xi_m^m(t)
\bigr), t\geq0\bigr).
\end{equation}
This fact is a simple consequence of the (lookdown) updating rule (\ref
{Epushup}) that makes
the type at level $i$ depend only on the previous types at levels up to
(and including)~$i$.
Therefore, one can unambiguously define
the label process $(\xi_i, i=1,2,\ldots)$ simultaneously
for all $i$, as
\begin{equation}
\label{defxi} \xi_i(t):= \xi_i^i(t) \equiv
\lim_{n\to\infty} \xi_i^n(t)\qquad  \forall t
\geq0,  \forall i\geq1.
\end{equation}

\begin{defn}\label{Dlookdownlabel}
We call $\xi:=(\xi_{i}(t),t\ge0)_{i\ge1}$ the \emph{label process}
associated to~$\pi$.
We may write $\xi^{\pi}$ for $\xi$ in order to indicate this association.
Unless otherwise specified we always assume that the $(\xi_i(0))_{i
\ge1}$ are i.i.d. uniformly distributed on $[0,1]$.
\end{defn}

In the sequel we will often focus on $(N^\pi(t), t \geq0)$, the
number of (distinct) types in the population process, defined by
\begin{equation}
\label{NX} N^\pi(t):=\#\bigl\{\xi_1(t),
\xi_2(t),\ldots\bigr\},\qquad t\geq0.
\end{equation}
Note that $N^\pi(t) \in\{1, 2, \ldots\} \cup\{\infty\}$ and $N^\pi
(0) = \infty$, due to our assumptions on $\xi(0)$.

The next proposition justifies the above definition of $\xi$ and
ensures that the corresponding limiting empirical measure exists
%jb the process is not always Markov so I changed the next sentence.
(as a c\`adl\`ag Markov process when the process $\pi$ is a Poisson
point process). These facts will be used in the construction of the
coupling without further reference in the sequel.

\begin{proposition}\label{Plookdowngeneral}
Let $\pi$ be a point process satisfying \eqref{EconditionH}, and
let $(\xi_i)_{i\ge1}$ be its label process.
Then the limit $\Xi_t=\lim_{n\to\infty} \frac{1}n \sum_{i=1}^n
\delta_{\xi_i(t)}$ exists simultaneously for all $t$ almost surely
and is c\`adl\`ag with respect to the weak topology.

Moreover, if $\pi$ is random and satisfies \eqref{EconditionH}
almost surely, $(\Xi_t,   t\geq0)$ is a Markov process in its own
filtration provided $U(t) = \sum_{i: t_i\le t} p_i^2$ has independent
increments.
\end{proposition}

\begin{defn}\label{Dlookdown}
The process $(\Xi_t,t\ge0)$ is the \textit{lookdown \textup{(}measure-valued\textup{)}
process} associated to $\pi$.
We may write $\Xi^{\pi}_t$ instead of $\Xi_t$ to make explicit the
dependence on the point process $\pi$.
\end{defn}

\begin{pf*}{Proof of Proposition~\ref{Plookdowngeneral}}
%jb changed wording. Should we signal in the proposition that this is
%done in dk99 ?
The proof can essentially be found in \cite{dk99}, up to a few
modifications due to the difference in points of view.
%in a more general setting which does not make clear the role of the
%point process $\pi$ and condition$(\ref{EconditionH})$ and
We explain how to adapt their arguments to our setting.
Recall the notation of Definition~\ref{Dlookdownlabel}.
To show that $\xi^n$ is well defined, note that, almost surely,
\begin{equation}
\label{Efinmany} \#\bigl\{j \geq1\dvtx t_j\in[0,t] \mbox{ and }
m_j(n)\ge2\bigr\} < \infty\qquad \forall t\geq0.
\end{equation}
Indeed, for each $j$ the indicator $\indica{m_j(n) \ge2}$ has expectation
$
1-(1-p_j)^n -n p_j(1-p_j)^{n-1} \leq{n \choose2} p_j^2
$,
and assumption $(\ref{EconditionH})$ together with Borel--Cantelli
lemma ensures~(\ref{Efinmany}).
Thus the dynamic (inductive) update (\ref{Epushup}) is feasible, and
the label process $\xi$ associated to $\pi$ is well defined.
A crucial feature of $\xi$
is that for each fixed $t>0$, the sequence
$(\xi_i(t),i=1,2,\ldots)$ is exchangeable.
Indeed, $(\xi_i(0),i=1,2,\ldots)$ is an exchangeable family, and the
transitions
preserve the exchangeability.
An application of de Finetti's theorem now yields the existence of the limit
\begin{equation}
\label{defXi} \Xi_t=\lim_{n\to\infty}
\Xi^n_t\qquad\mbox{where } \Xi^n_t:=
\frac{1}n \sum_{i=1}^n
\delta_{\xi_i(t)}
\end{equation}
for any fixed time $t$, and hence for all $t\in\Q$ simultaneously,
almost surely.

To see that the limit $\Xi_t$ actually exists simultaneously for all
$t$ with probability one is more delicate and is proved by Donnelly and
Kurtz in \cite{dk99}. Essentially one can adapt the proof of their
Lemma~3.4 to see that for each fixed $T>0, \varepsilon>0$ and each Borel
bounded function $f \dvtx [0,1] \mapsto\R$, there exists a positive
sequence $(\delta_l)_{l >0}$ such that $\sum_{l\ge0} \delta_l
<\infty$ and such that for all $l,m \ge k$
\begin{equation}
\P \biggl( \sup_{t\le T} \biggl\llvert \int_0^1
f(x) \Xi^m_t(dx) - \int_0^1
f(x) \Xi^l_t(dx) \biggr\rrvert >\varepsilon \biggr) \le
\delta_k.
\end{equation}
This implies that the sequence $\int_0^1 f(x) \Xi^m_t(dx)$ is almost
surely Cauchy.
Since the space of bounded measurable functions is separable (see,
e.g., Lemma~1.2 in \cite{dk99}) this is enough to guarantee existence
$(\Xi_t,t\ge0)$ as a process with values in the set of Borel
measures. Moreover $\Xi^n_t$ converges for all $t\le T$ simultaneously
almost surely.

Now assume that $\pi$ is a random point process satisfying \eqref
{EconditionH} and that $(U(t), t \ge0)$ has independent increments.
Then it is easy to check that the label process $(\xi_i(t), t\ge0)$
is Markov (in its own filtration). The Markov property for $\Xi$ then
follows directly from the fact that exchangeable laws on $[0,1]^\infty
$ are by De Finetti's theorem in one-to-one correspondence with the law
of their empirical measure on $[0,1]$. [Note that, however, $(\Xi_t,
t\ge0)$ is \textit{not} Markov with respect to the strictly greater
filtration of the label process, since the type of individual 1 will
tend to take over the population as time evolves.] A similar argument
is used in \cite{foucart}, Proposition~3, to prove the Feller property
for generalized Fleming--Viot processes with immigration. The context
there is slightly more general since the case of the so-called $\Xi
$-Fleming--Viot where simultaneous resampling events are allowed is
considered. However, the driving point measure used in \cite{foucart}
is Poissonian whereas we authorize more general point processes.
\end{pf*}

\begin{remark}
Donnelly and Kurtz \cite{dk99} work under a different set of
assumptions. Their setup is more general in the sense that they do not
assume the consistency of the finite-$n$ label processes $(\xi^n_i(t),
t \ge0)_{1 \le i \le n}$. (Furthermore, note that they also include a
Markov mutation diffusion operator that drives the motion of labels in
between reproduction events.) In fact, the total number of particles is
allowed to vary in their setting.
For this reason, their construction does not make sense conditionally
given the (limiting) point process $\pi$, which is an important
feature of our construction. The main novelty in our setting is the
observation that the assumption \eqref{EconditionH} is in fact all
that is needed to guarantee existence of the measure-valued process
$(\Xi_t, t \ge0)$ (it is also clear that this condition is necessary
for the very construction of the label process). In the notation of
Donnelly and Kurtz, this amounts to checking that the process $(U^n(t),
t \ge0)$ converges in distribution to $(U(t), t \ge0)$.
\end{remark}

\subsection{\texorpdfstring{Ancestral partitions, Fleming--Viot processes and $\Lambda$-coalescents}
{Ancestral partitions, Fleming--Viot processes and Lambda-coalescents}}\label{SFleVio}

We next apply Proposition~\ref{Plookdowngeneral}
in two different settings, corresponding to
the Fleming--Viot process and to the CSBP, respectively.
The upshot of this construction is a convenient way to track the
respective genealogies.
This is achieved through the \textit{ancestral partition process},
associated to the process $\xi$ constructed in Proposition~\ref
{Plookdowngeneral}.

Let $\pi$ be a point process satisfying \eqref{EconditionH}, and
$\xi^\pi$ its associated label process.
Note that for each $s > 0$, the shifted point process $\pi^{-s}:= \{
(p,t-s)\dvtx (p,t)\in\pi,   t\geq s\}$ also satisfies
\eqref{EconditionH}, and that, due to the updating rule (\ref
{Epushup}),
the label updates of the associated label process $\{\xi^{\pi
^{-s}}(t),  t \geq0\}$
are the same as those of $\{\xi^\pi(t),  t \geq s\}$. The difference
between the two processes is manifested through their initial states,
since for $i\neq j$
we have $\xi_i^{\pi^{-s}}(0)\neq\xi_j^{\pi^{-s}}(0)$, almost surely,
while it is possible that $\xi_i^\pi(s)=\xi_j^\pi(s)$.
Now fix some $T>0$.

\begin{defn} \label{Dancestralpartition}
%Let $(\xi(t),t\ge0), (A_i(j))_{i,j \in\N}, (m_i(j))_{i,j \in\N} $
%be as above.
%and fix $T>0.$
The ancestral partition process $(\mathcal{R}^T(t), 0\le t \le T)$
takes values in the space of level partitions (or partitions of $\N$).
% the partition-valued process where
For each $t\le T$,
$\mathcal{R}^T(t)$ is defined by the equivalence relation: $i\sim j$
in $\mathcal{R}^T(t)$ if and only if $\xi_i(T)$ and $\xi_j(T)$
descend from the same level at time $t$, or equivalently, if $\xi
_i^{\pi^{-t}}(T-t) = \xi_j^{\pi^{-t}}(T-t)$; see also equation (2.3)
in \cite{7}.
\end{defn}
Note that $\mathcal{R}^T(T)$ is the trivial partition $\{\{i\}\dvtx i \in
\N\}$ and that
$\mathcal{R}^T(t_1)$ is a coarser partition than $\mathcal
{R}^T(t_2)$, whenever $0\leq t_1\leq t_2 \leq T$.

We now briefly recall the definition of generalized Fleming--Viot
processes as well as their link to $\Lambda$-coalescents.
A generalized $\Lambda$-Fleming--Viot
process (in the sense of Bertoin and Le Gall \cite{blg2}) $(\rho_t, t
\ge
0)$ is a Markov process taking values in the space
$\mathcal{M}$ of probability measures on $[0,1]$. Its
generator $L$ is defined as follows: given a finite measure $\Lambda$
on $[0,1]$,
\begin{equation}
LF(\mu)=\int_{(0,1]}y^{-2}\Lambda(dy)\int
_{[0,1]}\mu(dx) \bigl( F\bigl((1-y)\mu+y\delta_{x}
\bigr)-F(\mu) \bigr), \label{GenFV}
\end{equation}
where $F\dvtx \mathcal{M}\to\R$ is a bounded continuous function.
In words, a number $y$ between 0 and 1 is
sampled at rate $y^{-2}\La(dy)$. A type $x$ is sampled from
$\rho_{t-}$.
Then $\rho_t$ is obtained from $\rho_{t-}$ by
scaling down $\rho_{t-}$ by $(1-y)$ and adding to the result an atom
at $x$ of mass $y$.
%Hence $\rho_t$ stays a probability measure.

%nb I am removing the following paragraph, as I do not see why it is
%needed here.
%Moreover, Bertoin and Le Gall \cite{blg1} show that
%generalized Fleming--Viot processes are characterized by the
%following property of their generators. Let $p \in\N$ and let $f$
%be a continuous function on $[0,1]^p$. Consider the function $G_f$
%on $\mathcal{M} $ defined by
%$$
%G_f(\mu)=\int_{[0,1]^p}\mu(dx_1)\ldots\mu(dx_p)
%f(x_1,\ldots,x_p).
%$$
%Then the action of the generator $L$ on the function $G_f$ is
%LG_f(\mu) = \sum_{I \subset\{1,\ldots,p\}, |I|\ge2} &
%where $R_I(x_1, \ldots, x_p)=(y_1,\ldots, y_p)$ is defined by
%$$
%y_j=x_{\min I} \mbox{ when $j\in I$ and } y_j=x_j \mbox{
%otherwise,}
%$$
%and where the coefficients $\lambda_{p,|I|}$ are the coalescence rates
%associated to the $\Lambda$-coalescent, i.e., $\lambda_{b,k} =

\begin{theorem}
\label{TinforGFV}
Let $\Lambda$ be a finite measure on $[0,1]$. Let $\pi$ be a Poisson
point process on $[0,1]\times\R$
with intensity $x^{-2} \Lambda(dx) \otimes dt$.
Then the lookdown process $\Xi^\pi$ (cf. Definition~\ref
{Dlookdown}) is a $\La$-generalized Fleming--Viot process with generator
\eqref{GenFV}, started from the uniform measure on $[0,1]$.
Furthermore, the ancestral
partition process $(\cR^T(T-t), 0\le t \le T)$ is the $\Lambda
$-coalescent, run for time $T$.
\end{theorem}

\begin{pf}A careful proof of this fact can be found in Lemma~3.6 of
\cite{7}, that is directly based on
the work of Donnelly and Kurtz \cite{dk99}.
We include a simpler proof which relies instead on the duality
introduced by Bertoin and Le Gall \cite{blg1}. We start with the claim
that the ancestral partition
process $(\mathcal{R}^T(T-t),0 \le t\le T)$ is the $\Lambda$-coalescent.
This follows
simply from the following observation: let $\pi'$ be the
point process obtained from $\pi$ by applying the transformation
$(p,t) \mapsto(p,T-t)$. Then $\pi'$ has same law as $\pi$ restricted
to $[0,T]$
and is thus a Poisson point process on $[0,1]\times[0,T]$ with
intensity $x^{-2} \Lambda(dx) \otimes dt\indic{[0,T]}(t)$.
Now, the updating rule~(\ref{Epushup}) can be rephrased as follows:
at each atom $(x,t)$ of $\pi'$ one flips a coin for each active
ancestral lineage with probability of heads equal to $p$ and the
lineages that come
up heads merge. This is precisely the Poisson process construction of
$\Lambda$-coalescents; see, for example, Theorem~3.2 in \cite{ensaios}.

Let $(\rho_t,t\ge0)$ be a Fleming--Viot process, and let
\[
F_t(x)=\rho_t\bigl([0,x]\bigr),\qquad  0\le x\le1
\]
be the associated \textit{bridge} process. Denote by $F_t^{-1}$ the
c\`adl\`ag inverse of the map $x \mapsto F_t(x)$.
%We start by recalling the following fact.
Let $V_1,V_2,\ldots,$ be
i.i.d. uniform random variables in $[0,1]$, independent
of $(\rho_t,t\ge0)$.
By the Glivenko--Cantelli theorem (see, e.g., (7.4) in
Chapter~1 of Durrett \cite{durrett}),
noting that
$(F_t^{-1}(V_i),i\ge
1)$ are i.i.d. samples from the random measure $\rho_t$,
we have for each fixed $t\ge0$
\begin{equation}
\rho_t = \lim_{n\to\infty}\frac{1}n\sum
_{i= 1}^n \delta_{F_t^{-1}(V_i)}\qquad\mbox{almost surely,} \label{Erhothree}
\end{equation}
where the limit is taken in the sense of the weak topology on
probability measures.

Bertoin and Le Gall \cite{blg2} proved
that the $\Lambda$-coalescent $(\Pi_t,t\ge0)$ is dual to the
generalized
Fleming--Viot process corresponding to $\La$ in the following sense:
if $n\ge1$
and $f$ is any continuous function on $[0,1]^n$, then
\begin{equation}
\label{duality} \E\bigl(f\bigl(F^{-1}_t(V_1),
\ldots,F^{-1}_t(V_n)\bigr)\bigr) = \E\bigl(f
\bigl(Y\bigl(\Pi^n(t),V'_1,\ldots,
V'_n\bigr)\bigr)\bigr),
\end{equation}
where $\Pi^n(t)$ denotes the restriction of $\Pi_t$ to $[n]$,
the random variables $(V'_1,\ldots,\break V'_n)$ are i.i.d. uniform
on $[0,1]$, and independent of $(\Pi_t,t\ge0)$, and where the map $Y$ is
defined as follows:
\[
\begin{tabular}{p{270pt}@{}}
$\mbox{for }\pi\in\mathcal{P}_n\mbox{ and }(x_1,
\ldots,x_n)\in [0,1]^n \mbox{ let }
Y(\pi, x_1, \ldots, x_n) = (y_1, \ldots,
y_n) \mbox{ with } y_j = x_i \mbox{ for }
i=\min\{ k \dvtx k \sim_{\pi} j \}$.
\end{tabular}
\]
Note that the duality relation (\ref{duality}) has the form of a generalized
functional duality in the context of
interacting particle systems (see \cite{liggett}), and
should not be confused with the notion of duality between
coagulation and fragmentation processes of~\cite{saintflour}.

We next verify that, for each $t>0$,
\begin{equation}
\label{sample1} \bigl(\xi_1(t),\ldots, \xi_n(t)\bigr)
\eqlaw Y\bigl(\Pi^n(t),V_1',
\ldots,V_n'\bigr).
\end{equation}
This fact is an immediate consequence of our construction. Indeed, at
time $t$ two levels $i$ and $j$ have the same type $\xi_i(t)=\xi
_j(t)$ if and only if they descend from the same level at time 0 [since
all the $\xi_i(0)$ are almost surely distinct]. Hence $\xi_i(t)=\xi
_j(t)$ if and only if $i$ and $j$ belong to the same block of $\mathcal
{R}^t(0)$.
Therefore
\[
\bigl(\xi_1(t),\ldots, \xi_n(t)\bigr) = Y'
\bigl(\mathcal{R}^t(0),\xi_1(0),\ldots,
\xi_n(0)\bigr),
\]
where for $\pi=(B_1,B_2,\ldots) \in\mathcal{P}_n$
and $(x_1,\ldots,x_n)\in[0,1]^n$ we let
\begin{eqnarray*}
Y'(\pi, x_1, \ldots, x_n) =
(y_1, \ldots, y_n) \qquad\mbox{with } y_j =
x_i \mbox{ for } j\in B_i.
\end{eqnarray*}
Clearly, as long as the random variables $\Pi\in\P_n$ and
$(X_1,\ldots, X_n)$ (i.i.d. uniform on $[0,1]$) are independent one has
\[
Y(\Pi,X_1,\ldots,X_n) \eqlaw Y'(\Pi,
X_1, \ldots, X_n),
\]
and since the $\xi_i(0)$ are i.i.d. uniform on $[0,1]$ and $\mathcal
{R}^t(0) \eqlaw\Pi(t)$, this proves the claim (\ref{sample1}). Due to
(\ref{duality}), one concludes that
$(F^{-1}_t(V_1),\ldots,F^{-1}_t(V_n))_{t \ge0} $ and $(\xi_1(t),
\ldots,\xi_n(t))_{t \ge0}$ have the same one-dimensional
marginals.
This implies that
\[
\forall t\geq0\qquad \Xi^n_t \eqlaw\frac{1}n\sum
_{i=1}^n \delta_{F_t^{-1}(V_i)}, \qquad  n\geq1\quad
\mbox{and hence that}\quad \Xi _t\eqlaw\rho_t.
\]
Our argument was carried out under the assumption that the initial state
is the uniform law on $[0,1]$. However, it would equally apply if the
$\xi_i(0)$ were drawn independently from any other law on $[0,1]$.
Since $\Xi$ and $\rho$ are both c\`adl\`ag Markov processes,
they must be equal in distribution.
\end{pf}

\subsection{Lookdown process of a CSBP}
Recall $\psi$ from (\ref{Dpsi}) and
consider a CSBP $(Z(t),t\ge0)$ with branching mechanism
$\psi$; see, for example, \cite{athreya-ney} or \cite{ensaios}, Chapter~4.2, for an elementary introduction.
In the sequel, we often refer to any such process as $\psi$-CSBP.
In this section assume that $Z$ is started from $Z(0)=1$.
Following Bertoin and Le Gall \cite{blg0}, recall existence of a two
parameter branching family $(Z_t(x), t \ge0, x \in[0,1])$,
such that for each fixed $x\in[0,1]$, $(Z_t(x), t\geq0)$ is a $\psi
$-CSBP started from $Z_0(x) = x$, independent from the $\psi$-CSBP
$(Z_t(1) - Z_t(x),   t\geq0)$.
In particular $(Z_t(1), t\geq0)\eqlaw(Z(t), t \geq0)$.
The quantity $Z_t(x)$ can be interpreted as the population size at time
$t$, descended from the initial fraction $x$ of the population at time $0$.
Furthermore, the branching property also implies that, for any $t>0$,
$(Z_t(x), x\in[0,1])$ is a subordinator.

We briefly recall the setting of \cite{7}.
For each fixed $t\geq0$, define $M_t([x_1,\break x_2]):= Z_t(x_2) - Z_t(x_1)$,
for all $0\leq x_1\leq x_2\leq1$. Then $M_t$ extends to a random
measure on $[0,1]$.
The process $M=(M_t, t\geq0)$ is easily seen to be Markov, with a
generator given by (see (1.15) in \cite{7} for the general case formula)
\[
\mathcal{L}F(\mu) = \int_0^1 \mu(da) \int
_{[0, 1]} \nu(dh) \bigl(F(\mu +h\delta_a)-F(
\mu)-hF'(\mu;a)\bigr),
\]
where $\nu(dh) = \La(dh)/h^2$, and $F'(\mu;a)$ denotes the Fr\'echet
derivative of $F$ at $\mu$ in the direction $\delta_{a}$; see, for
example, (1.4) in \cite{7}.
% (the term involving $F'(\mu;\cdot)$ comes from the intrinsic drift
%compensator, in analogy to a similar term for certain {\Levy}
%processes).
The process $M$ encodes the genealogy of the CSBP $(Z_t(1), t\ge0)$
(this is a continuous time/space analogue to the relation between a
Galton--Watson process and the associated tree).
The composition of the population is then well described by the \textit{ratio process}
$R=(R_t, t\geq0)$ defined by $R_t= \frac{1}{Z_t(1)}M_t$, taking
values in the space of probability measures.
Now define
\begin{equation}
\label{DpiZ} \pi^Z= \bigl\{ \bigl({\Delta Z(t)}/{Z(t)},t \bigr)\dvtx t\geq0 \bigr\}
\end{equation}
to be the point process of normalized jump sizes of $Z$.
Here and below (without further mention), we will account in $\pi^Z$
only the points
$({\Delta Z(t)}/{Z(t)},t) \in(0,\infty) \times[0,\infty)$, which
represent the true jumps of the process.
\begin{lemma}\label{TinforGFV+CSBP}
The condition \eqref{EconditionH} holds for $\pi^Z$, and the
associated lookdown process $(\Theta_t(\cdot), t \ge0)$ is equal in law
to the ratio process $(R(t),t\ge0)$.
%zzz in 7, this is done only if \sigma=0 no Brownian part in psi. Est
%ce que on doit le prouver ?
%zzzvl we also assume that the BM part is 0 in this whole section
\end{lemma}

\begin{pf}
A detailed proof is given in the ``Proof of (2.4)'' in \cite{7}, pages
313--315, although the idea goes back at least to Theorem~3.2 in \cite{dk99}.
\end{pf}

\subsubsection{Evolution of the number of types}
Let $Z$ be a CSBP with branching mechanism $\psi$ started from
$Z_0=1$, and assume that Grey's condition (\ref{Egrey0}) is satisfied.
Denote by $\zeta_Z =\inf\{t \ge0 \dvtx Z(t)=0 \}$ its (almost surely
finite) extinction time.
Let $\pi^Z$ be the associated point process of rescaled jump sizes
\eqref{DpiZ}, and
note that $\pi^Z$ has no points in $[0,1]\times(\zeta_Z,\infty)$.
Recall definition (\ref{NX}),
and define $N^Z(t) = N^{\pi^Z}(t)$, $t<\zeta_Z$, and $N^Z(t)=0$,
$t\ge\zeta_Z$.

%jb put definition of v here with citation of duleg
Let us define $v(t):= \inf\{z> 0\dvtx \int_z^\infty\psi(q)^{-1}\,dq <t\}
$ with the convention that $\inf\varnothing=\infty$ or equivalently
let $v(t)$ be the solution of
\begin{equation}
\label{Dv} \int_{v(t)}^\infty\frac{dq}{\psi(q)}
=t.
\end{equation}
Recall from Duquesne and Le Gall \cite{duleg} that the function $v$
describes the evolution of the number of alive families at time $t$ in
a $\psi$-CSBP. More precisely, we have the following:
\begin{proposition}
\label{PNXandNZ}
If (\ref{Egrey0}) is satisfied, then $N^Z(t)< \infty$, for all
$t>0$, almost surely, and moreover
\begin{equation}
\bigl(N^Z(t),t\ge0\bigr) \eqlaw\bigl(Q \bigl(v(t)\bigr), t\ge0
\bigr),
\end{equation}
where $t \mapsto Q(t)$ is a standard Poisson counting process, and
where $v(t)$ is defined in (\ref{Dv}).
In particular,
\begin{equation}
\lim_{t\to0}\frac{N^Z(t)}{v(t)} = 1\qquad \mbox{ almost surely.}
\end{equation}
If (\ref{Egrey0}) is not satisfied, then both $v$ and $N^Z$ are
infinite for all $t>0$ almost surely.\vadjust{\goodbreak}
%jb added discussion of case v=\infty
\end{proposition}
\begin{pf}
This essentially follows from
Theorem~12 in \cite{bbs2} and Corollary~1.4.2(ii) in \cite{duleg};
see also Corollary~4.1 in \cite{ensaios} for an elementary sketch of proof.
Indeed, when Grey's condition is satisfied, we may
use the construction of \cite{bbs2} for the Donnelly--Kurtz
lookdown process, where the labeling process
$(\theta_1(t),\theta_2(t),\ldots)$ is directly defined in terms of
the excursions of a Continuous Random Tree (CRT) with branching
mechanism $\psi$; see \cite{duleg} or \cite{bbs2} for the basic
terminology and
properties of these objects, to which we will refer in this
proof. Let $(H_s,0\le s\le T_1)$ be the height process associated
with $(Z_s,s\ge0)$, where $T_1:= \inf\{u>0\dvtx L_u^0>1\}$ and
where $(L_u^0,u \ge0)$ is the local time process at level 0 of
$(H_s,s\ge0)$. It follows from
the construction in \cite{bbs2} that one can embed the lookdown
construction in the CRT so that for any $t>0$,
$N^Z(t)$ is exactly the number of
excursions of $(H_s,0\le s\le T_1)$ that reach level $t$. It follows
directly [by excursion theory for $(H_s,0\le s \le T_1)$] that
$(N_Z(t),t>0)$ has the law
of $(Q_{\tilde v(t)},t\ge0)$, where by definition,
\[
\tilde v(t) = \mathbf{N}\Bigl(\sup_{s\ge0} H_s >t
\Bigr).
\]
Here, $\mathbf{N}(\cdot)$
denotes the excursion measure of $H$. By Corollary~1.4.2(ii) of \cite{duleg},
$\tilde v(t) = v(t) < \infty$, which proves the result.
\end{pf}

\begin{remark}
For each fixed $t>0$, due to the exchangeability of\break  the~sequence $(\xi
_i(t), i=1,2, \ldots)$, the number of types $N^Z(t)$ is almost surely
equal to the~number of atoms of the
purely atomic measure $\Xi_t =\break \lim_{n\to\infty} n^{-1} \sum_{i=1}^n \delta_{\xi_i(t)}$.
\end{remark}

\begin{remark}
The property $\P(N^Z(t)<\infty)=1$ may seem counter-intuitive in view
of the fact that types are not destroyed in any particular application
of the updating rule (\ref{Epushup}).
However, an accumulation of many densely placed small lookdown jumps
``pushes off'' to infinity all but finitely many types in any positive
amount of time,
whenever Grey's condition \eqref{Egrey0} is fulfilled.
\end{remark}

\section{The coupling}
\label{Scoupling}

\subsection{Coupling construction}
\label{SScouplingconstr}

We can now explain the coupling between $\Lambda$-coalescents and
CSBP. The key idea is to use the following result due to Lamperti,
which expresses
any CSBP as a time-change of a L\'evy process.

Consider a L\'evy process $(X_t, t\ge0)$ with Laplace exponent $\psi
$ given in (\ref{Dpsi}), and assume $X_0 =x\in(0,1]$.
Define
\begin{equation}
\label{Elamperti2} U^{-1}(t):= \inf \biggl\{ s>0 \dvtx \int
_0^s \frac{du}{X_u} > t \biggr\}
\end{equation}
and
\begin{equation}
\label{Elamperti1} Z_t = X_{U^{-1}(t)},\qquad  t\ge0.
\end{equation}

\begin{theorem}[(Lamperti \cite{lamp,lamperti2})]
The process $(Z_t,  t\ge0)$ is a $\psi$-CSBP started from $Z_0 = x$.
\end{theorem}

\textit{Construction}.
We now describe the coupling between the genealogies of a CSBP and
Fleming--Viot processes.
Assume that the L\'evy process $X$ and its corresponding CSBP $Z$
(Lamperti time-changed as above) satisfy $X_0=Z_0=1$.
As before, denote by $\pi^Z$ the point process of the rescaled jump
sizes of $Z$.
Call $\xi= (\xi_i(t), t\ge0)_{i\ge1}$ the label process of $\pi^Z$
obtained from the lookdown construction applied to $Z$.

Consider simultaneously the point process $\pi^X = (\Delta X(t_i),
t_i)$ of (unscaled) jump sizes of $X$, and its
associated label process $\theta= (\theta_i(t), t\ge0)_{i\ge1}$,
as well as the lookdown measure $\Theta=(\Theta_t, t\ge0)$. Then
$\Theta$ is a $\Lambda$-Fleming--Viot process, and hence
(due to Theorem~\ref{TinforGFV}) has a genealogy given by a
$\Lambda$-coalescent.
Indeed, since $X$ is a L\'evy process, due to the L\'evy--It\^o decomposition,
the point process of jumps $\pi^X = (\Delta X({t_i}), t_i)$ is a
Poisson point process with intensity $\nu(dx) \otimes dt$, where $\nu
(dx)=x^{-2}\Lambda(dx)$ is the
L\'evy measure of $X$.

\textit{Heuristics}. For a small $t>0$,
the two point processes $\pi^X$
and $\pi^Z$, restricted to $[0,t]$, are ``close to each other.''
Indeed, each point
$(p,t) \in\pi^X$ also corresponds to a point $(\tilde p, \tilde t)
\in\pi^Z$, where $ t = U^{-1}(\tilde t)$, and $\tilde p = p/Z(\tilde
t)$. Now, since $(X_t,t\ge0)$ is almost surely continuous at $t=0$,
the time-change $U^{-1}$ is almost surely differentiable at $t=0$ with
derivative close to $1$. Therefore, $U^{-1}(t) \sim t$ as
$t \to0$, and one deduces that for small $t$, $\tilde t \approx t$.
Likewise, invoking the continuity of $Z$ and the fact $Z_0=1$, we have
$Z(\tilde t) \approx1$, hence
$(\tilde p, \tilde t) \approx(p,t)$.

It is therefore reasonable to believe that for small $t$,
$N^X(t) \approx N^Z(t)$, where $N^X(t)$ [resp.$,  N^Z(t)$] is
the number of types in the lookdown process associated to $\pi^X$
(resp., $\pi^Z$) at time $t$.
At the same time, by Proposition~\ref{PNXandNZ} we also know
$N^Z(t)\sim v(t)$ almost surely as $t\to0$, and all of the above
strongly suggests that the same is true for $N^X$ in place of $N^Z$.

Finally, due to Theorem~\ref{TinforGFV}, we have
\begin{equation}
\label{Eeqivlaw} N^X(t)\eqlaw N^\La(t)\qquad \mbox{for each
fixed $t \ge0$,}
\end{equation}
where $N^\La(t)$ is (as usual) the number of blocks in the
corresponding $\Lambda$-coalescent at time $t$. The reader can easily
check this property by restricting attention to the first $n$ levels,
and using the updating rule (\ref{Epushup}), as well as the fact
that $(\pi^X(t), t\in[0,T])$ and
$(\pi^X(T-t), t\in[0,T])$ have the same distribution.
Therefore, we obtain $N^\La(t)\sim v(t)$ in probability, as $t\to0$.

 We will now turn these heuristic observations into a rigorous
argument for Proposition~\ref{Tsmalltime},
starting with a monotonicity lemma.
%As indicated in the Introduction, this convergence in fact holds
%almost surely (see Theorem~1 of \cite{bbl1}).
%vl we already said this in the intro and we are saying it in Remark~16
%below

\begin{defn}\label{DorderPPP}
Given two point processes $\pi$ and $\pi^+$ on $[0,1]\times\R_+$ on
the same probability space,
and a random time $T \ge0$, measurable with respect to the filtration
generated by $\pi$ and $\pi^+$,\vadjust{\goodbreak}
we write $\pi\lhd_{|[0,T]} \pi^+$ (or~$\pi\lhd\pi^+$ on $[0,T]$)
if there exists an
increasing c\`adl\`ag process $r \dvtx [0,T] \mapsto\R^+$
such that, almost surely, $r(0)=0$ and
\[
\pi= \bigl\{(p_i,t_i)\dvtx i\geq1\bigr\} \quad\mbox{and}\quad \pi^+=
\bigl\{\bigl(q_i,r(t_i)\bigr)\dvtx i\geq1\bigr\},
\]
where $p_i\le q_i$, for each $i\ge1$ such that $t_i\leq T$.
\end{defn}
In words, $\pi\lhd\pi^+$ on $[0,T]$, if the atoms of $\pi^+$ are
those of $\pi$, time-changed by $r$ and multiplied in size by a
(possibly nonconstant and random) quantity not smaller than $1$.
Observe that $r$ preserves the order of the atoms, almost surely.
In our main applications, the form of $r$ will be rather simple.
Furthermore, the processes $\pi$ and $\pi^+$ of interest will both
have (countably) infinitely many atoms in any interval of positive
length, almost surely, ensuring that $\{r(T)<\infty\}=\{T<\infty\}$,
almost surely.

Consider now $\pi$ and $\pi^+$ such that $\pi\lhd\pi^+$ on $[0,T]$
for some finite random time $T$, and both
\[
\sum_{i: t_i\le t} p_i^2 <\infty,\qquad
\sum_{i: t_i\le t} q_i^2 <\infty \qquad \forall t\geq0, \mbox{ almost surely}.
\]
One can then construct a coupling of $\Xi^{\pi}$ [with its label
processes $\xi= (\xi_i(t),t\ge0)_{i\ge1}$] and $\Xi^{\pi^+}$
[with its label processes $\xi^+ = (\xi^+_i(t),t\ge0)_{i\ge1}$], by
using the same collection $\{U_{i,j}\}_{i,j\in\N}$ of i.i.d. uniform
random variables to specify the levels participating in the resampling
events in Definition~\ref{Dlookdownlabel}.
Due to $\pi\lhd\pi^+$ on $[0,T]$,
the following result is obvious by construction:

\begin{lemma}\label{monotonicity}
If $\pi\lhd\pi^+$ on $[0,T]$, then
\[
\P\bigl(N^{\pi^+}\bigl(r(s)\bigr) \le N^\pi(s)\ \forall s
\in[0,T]\bigr)=1.
\]
\end{lemma}

\subsection{\texorpdfstring{Proof of Theorem \protect\ref{Tmainthm} and the asymptotics for the number of blocks}
{Proof of Theorem 1 and the asymptotics for the number of blocks}}
\label{SappltoN}

To prove Theorem~\ref{Tmainthm} it suffices to show that $N^\La(t)$
is infinite for all $t>0$ whenever $v(t)=\infty, \forall t >0$ and is
finite for all $t>0$ in the converse case.
This is now a consequence of the above coupling, used to show the
following proposition.

%jb the Tsmalltime is now mainly seen as a tool to prove Thm 1 ans is
%now a proposition. Parts of the discussion which was in the intro has
%been moved here after the prop.
%
\begin{proposition} \label{Tsmalltime}
For each $\eps\in(0,1)$,
\begin{eqnarray}
\label{comp} &&\P \biggl( \liminf_{t \to0} \frac{N^X(t)}{v ({(1+\eps)}/{(1-\eps)} t )} \ge
\frac{1}{1+\eps} ,
\nonumber
\\[-8pt]
\\[-8pt]
\nonumber
&&\hspace*{-8pt}\qquad  \limsup_{t\to0} \frac{N^X(t)}{v ({(1-\eps)}/{(1+\eps)} t )} \le
\frac{1}{1-\eps} \biggr) =1,
\end{eqnarray}
and therefore
\begin{equation}
\label{EapplidecompaNLA}\quad  \lim_{t \to0} \P \biggl[ \frac{1}{(1+\eps)^2}\cdot{v
\biggl(\frac
{1+\eps}{1-\eps} t \biggr)}\le N^\La(t) \le
\frac{1}{(1-\eps)^2}\cdot{v \biggl(\frac{1-\eps}{1+\eps} t \biggr)} \biggr] =1.
\end{equation}
\end{proposition}
\begin{remark}
\label{Rlimtech}
Observe that $N^X$ and $N^\La$ have only the same one-dimen\-sional
marginal distributions, but they are not equal in distribution as
processes. For instance while the first one only decreases by jumps of
size $1$ (this is known at least in the stable case; see \cite
{labbe}), the second one can decrease by jumps of arbitrary integral
length. Thus one cannot obtain more than (\ref{EapplidecompaNLA})
from (\ref{comp}).
This result is clearly weaker than Theorem~1 in \cite{bbl1},
\begin{equation}
\label{Talmostsure} \lim_{t \to 0} \frac{N^\La(t)}{v(t)}= 1 \qquad\mbox{almost surely}.
\end{equation}
As mentioned in the \hyperref[sec1]{Introduction}, it is the use of a sophisticated
martingale technique which yields this stronger result there.
However, it was the knowledge of the coupling described below that
initiated \cite{bbl1} and suggested the
form of the asymptotics in the first place.
\end{remark}

\begin{pf*}{Proof of Proposition~\ref{Tsmalltime}}
We start by showing (\ref{comp}) for $\eps$ sufficiently small. The
conclusion (\ref{EapplidecompaNLA}) will then readily follow.
Let us assume for the moment that ${\rm supp}(\Lambda)\subset[0,\eta
]$ where $\eta<1$,
and fix some $\eps\in(0,1/\eta-1)$.
Consider again the L\'evy process $X$ with Laplace exponent $\psi$
such that $X_0=1$, and let
\[
\pi= \pi^X = \bigl\{(\Delta X_t,t)\dvtx t >0\bigr\}
\]
be the corresponding Poisson point process. Let $\pi_\eps^-$
(resp., $\pi_\eps^+$) be the image of $\pi$ under the
map $(p,t) \mapsto( p(1-\eps),t)$ [resp., $(p,t) \mapsto
(p(1+\eps),t)$].
Due to our assumptions on ${\rm supp}(\Lambda)$ and the choice of
$\eps$,
we have that for each atom $(p,t)$ of $ \pi$,
$p(1+\eps)<1$ almost surely.
Therefore, both $\pi_\eps^+$ and $\pi_\eps^-$ are Poisson point
processes on $ (0,1) \times\R_+$.
Let $ \nu_{\eps}^+ \otimes dt$ (resp., $\nu_\eps^- \otimes dt$) be
the intensity
measure corresponding to $\pi_\eps^+$ (resp., $\pi_\eps^-$). If $f$
is a Borel function on $[0,1]$, then $\nu^+_\eps$
is obtained by the formula
\[
\int_{[0,1]} f(x)\nu_\eps^+(dx)=\int
_{[0,1]} f\bigl(x(1+\eps)\bigr)\nu(dx),
\]
and $\nu_\eps^-$ is obtained by an analogous formula with $1-\eps$
in place of $1+\eps$. For
$\lambda>0$, let
\[
\psi^\pm_\eps(\lambda):=\int_{(0,1)}
\bigl(e^{-\lambda x}-1+\lambda x\bigr) \nu^\pm_\eps(dx).
\]
By the above observation we see that, for each $\lambda>0$,
\begin{equation}
\label{malin} \psi^+_\eps(\lambda)=\psi\bigl(\lambda(1+\eps)\bigr),\qquad
 \psi^-_\eps (\lambda)=\psi\bigl(\lambda(1-\eps)\bigr).
\end{equation}
Therefore, if we let
$u_\eps^\pm(t):=\int_t^\infty \,d\lambda/\psi^\pm_\eps(\lambda)$
and
%jb v_\eps is the cadlag inverse of u_\eps, not just 1/u no ?
$v_\eps^\pm(t)= (u_\eps^\pm)^{-1}(t)$ the c\'{e}dl\'{e}g inverse of
$u_\eps^\pm$, we
have
\[
u^+_\eps(s)=\frac{1}{1+\eps}u\bigl(s(1+\eps)\bigr) \quad\mbox{and}\quad
u^-_\eps(s)=\frac{1}{1-\eps}u\bigl(s(1-\eps)\bigr),
\]
hence
\begin{equation}
\label{veps} v^+_\eps(t)=\frac{1}{1+\eps}v\bigl(t(1+\eps)\bigr)
\quad\mbox{and}\quad v^-_\eps(t)=\frac{1}{1-\eps}v\bigl(t(1-\eps)\bigr).
\end{equation}
Recall that $X_0=1$, and define
\[
X^+_t = (1+\eps) X_t -\eps \quad\mbox{and}\quad
X^-_t = (1-\eps) X_t +\eps \qquad  t>0.
\]
Then it is easy to see that both $(X^+_t,t\ge0)$ and $(X^-_t,t\ge0)$
are L\'evy processes such that $X^+_0= X^-_0=1$.
Moreover, the Laplace exponent of $X^+$ (resp., $X^-$) is
$\psi^+_{\eps}$ (resp., $\psi^-_{\eps}$).

Define
$T_\eps^+=\inf\{ s \dvtx |X^+(s)-1|> \eps\}$ and
$T_\eps^-=\inf\{ s \dvtx |X^-(s)-1| >\eps\}$.
Then, almost surely we have, for all $t\geq0$
\begin{eqnarray}
\label{comp+-} \frac{\Delta X^-(t)}{X^-(t)} \le \frac{\Delta X^-(t)}{1-\eps} = \Delta X(t) =
\frac{\Delta X^+(t)}{1+\eps} \le \frac{\Delta X^+(t)}{X^+(t)}
\nonumber
\\[-8pt]
\\[-8pt]
\eqntext{\mbox{on } \bigl\{t \le
T_\eps^+ \wedge T_\eps^-\bigr\}.}
\end{eqnarray}
Using the Lamperti transform, now define two continuous-state
branching processes with branching mechanism $\psi^+_\eps$ and
$\psi^-_\eps$, respectively, by setting $U_\pm(t):= \int_0^t \frac{1}{X^\pm_u}\,du$,
\begin{eqnarray}
U_\pm^{-1}(t):= \inf \bigl\{ s \ge0 \dvtx U_\pm(s) >t \bigr\}\quad \mbox{and}\quad Z^+_t:=
X^+_{{U_+^{-1}(t)}}, \qquad  Z^-_t:= X^-_{{U_-^{-1}(t)}},\nonumber\\
  \eqntext{t\geq0.}
\end{eqnarray}
Finally
%jb suppressed (compare with Lemma~\ref{TinforGFV+CSBP})
define $\pi^{Z^+}:= \{(\Delta Z^+_s/Z^+_s,s)\dvtx s\geq0\}$
and $\pi^{Z^-}:= \{(\Delta Z^-_s/Z^-_s,s)\dvtx\break  s\geq0\}$.
Due to (\ref{comp+-}), we have that almost surely
\begin{eqnarray*}
&&\pi^{Z^-} \lhd_{|[0,U_-(T_\eps^+\wedge T_\eps^-)]} \pi \qquad\bigl(\mbox {with
}r=U_-^{-1}\bigr) \quad\mbox{and}\\
 &&\pi\lhd_{|[0,T_\eps^+\wedge T_\eps^-]}
\pi^{Z^+} \qquad(\mbox{with }r=U_+),
\end{eqnarray*}
where $\lhd$ is as in Definition~\ref{DorderPPP}.
Both $T_\eps^+\wedge T_\eps^-$ and $U_-(T_\eps^+\wedge T_\eps^-)$
are clearly strictly positive and finite, almost surely.
Hence, Lemma~\ref{monotonicity} gives that almost surely, for all
$t\geq0$,
\[
N^\pi(t) \le N^{ \pi^{Z^-}}\bigl(U_-(t)\bigr)\quad \mbox{and}\quad
N^\pi(t) \ge N^{ \pi^{Z^+}}\bigl(U_+(t)\bigr)\qquad \mbox{ on } \bigl
\{t \leq T_\eps^+\wedge T_\eps^-\bigr\}.
\]
%
%jb introduced the following observation
Observe that this is already enough to prove Theorem~\ref{Tmainthm}
since $v_\eps^\pm$ is finite if and only if $v$ is finite, and thus
$N^{ \pi^{Z^+}}(U_+(t))=\infty$ for all $t>0$ if $v(t)=\infty$ for
all $t>0$ and likewise $N^{ \pi^{Z^-}}(U_-(t))<\infty$ for all $t>0$
if $v(t)<\infty$ for all $t>0$.

Proposition~\ref{PNXandNZ} implies that
\[
\lim_{t\to0}\frac{ N^{ \pi^{Z^-}}(t)}{v_\eps^-(t)} = \lim_{t\to0}
\frac{ N^{ \pi^{Z^+}}(t) }{v_\eps^+(t)} = 1 \qquad\mbox{almost surely.}
\]
This together with
$\P(T_\eps^+\wedge T_\eps^- >0)=1$ and the discussion above
yields
\begin{equation}
\limsup_{t\to0} \frac{N^\pi(t)}{v_\eps^-(U_-(t))} \leq1\quad \mbox{and}\quad \liminf
_{t\to0} \frac{N^\pi(t)}{v_\eps^+(U_+(t))} \geq1, \label{comp3-}
\end{equation}
almost surely.
Moreover, it is easy to check that almost surely, for all $t\geq0$,
\begin{equation}
\label{comp2+} t/(1+\eps)\le U_\pm(t)\le t/(1-\eps)  \qquad\mbox{on }
\bigl\{t \le T^+_\eps\wedge T^-_\eps\bigr\}.
\end{equation}
Due to monotonicity of $v_\eps^\pm$ and (\ref{comp2+}), we have that
again almost surely, for all $t\geq0$,
\begin{eqnarray}
\label{EMNhelp} \frac{N^\pi(t)}{v_\eps^-(U_-(t))} \ge \frac{N^\pi
(t)}{v_\eps^-(t/(1+\eps))}\quad \mbox{and}\quad
\frac{N^\pi(t)}{v_\eps^+(U_+(t))} \le\frac{N^\pi(t)}{v_\eps
^+(t/(1-\eps))}
\nonumber
\\[-8pt]
\\[-8pt]
\eqntext{\mbox{on } \bigl\{t \le
T^+_\eps\wedge T^-_\eps\bigr\}.}
\end{eqnarray}
Combining (\ref{veps}), (\ref{comp3-}) and (\ref{EMNhelp}), and
recalling $\P(T_\eps^+\wedge T_\eps^- >0)=1$,
we can now conclude that
\begin{eqnarray}
\label{Elasti} \limsup_{t \to0} \frac{N^\pi(t)}{v (t {(1-\eps)}/{(1+\eps)}
) } &\le&
\frac{1}{1-\eps} \quad\mbox{and}
\nonumber
\\[-8pt]
\\[-8pt]
\nonumber
\liminf_{t \to0}
\frac{N^\pi(t)}{v (t {(1+\eps)}/{(1-\eps)}
)}& \ge&\frac{1}{1+\eps}\qquad\mbox{almost surely.}
\end{eqnarray}
Since $N^X =N^\pi$ by definition, this gives (\ref{comp}), under the
hypothesis that
$\La$ does not give positive mass to a neighborhood of $1$.
Otherwise, we modify the above argument in the following way.
For a fixed $\eta\in(0,1)$, since $x^{-2} \La(dx)$ assigns a finite
mass to $(1-\eta,1]$, the first time $T_\eta$ when $X$ makes a jump
of size strictly greater than $\eta$ has an exponential random
variable law (with finite rate), hence it is strictly positive with
probability one.
The analysis (\ref{comp+-})--(\ref{Elasti}) clearly works if
$T^+_\varepsilon\wedge T^-_\varepsilon$ is everywhere replaced by
$T^+_\varepsilon\wedge T^-_\varepsilon\wedge T_\eta$,
yielding (\ref{comp}).

In particular, almost surely, for all $t$ sufficiently small,
\[
\frac{1}{(1+\eps)^2} \cdot v \biggl(t \frac{1+\eps}{1-\eps} \biggr) \le
N^{X}(t) \le\frac{1}{(1-\eps)^2} \cdot v \biggl(t \frac{1-\eps
}{1+\eps}
\biggr).
\]
The limit (\ref{EapplidecompaNLA}) is easily deduced from (\ref
{Eeqivlaw}) and this final estimate.
\end{pf*}

The asymptotics (\ref{Talmostsure}) in the sense of convergence in
probability can be obtained from Proposition~\ref{Tsmalltime} under additional
assumptions on $v$ (i.e., on $\La$) as the following result shows.
\begin{proposition}
\label{Palmostsure}
Assume $\La(\{0\})=0$. Then the convergence
\[
N^\La(t) /v(t) \to1 \qquad \mbox{in probability}
\]
holds at least if
\begin{equation}
\label{cond1} \lim_{\eps\to0} \limsup_{t\to0}
\frac{v(t(1-\eps))}{v(t)}=1,\qquad  \lim_{\eps\to0} \liminf_{t\to0}
\frac{v(t(1+\eps))}{v(t)}=1,\vadjust{\goodbreak}
\end{equation}
and, in particular, if
\begin{equation}
\label{cond2} \psi\bigl(v(t)\bigr) = O\bigl(v(t)/t\bigr)\qquad\mbox{as } t \to0.
\end{equation}
\end{proposition}

\begin{pf}%{Proof of Proposition~\ref{Palmostsure}}
The first claim follows by simple calculus manipulations from (\ref
{EapplidecompaNLA}).
To see why (\ref{cond2}) implies (\ref{cond1}),
we note that $\psi\dvtx [0,\infty) \to\R^+$ of (\ref{Dpsi}) is
a (strictly) increasing and convex function on $[0,\infty)$.
Furthermore,
$v_\psi'(s)=-\psi(v_\psi(s))$,
so that $v_\psi$ is decreasing with its derivative decreasing in
absolute value. Therefore, for $\eps>0$ small enough,
\[
\bigl|v\bigl(t(1+\eps)\bigr) - v(t)\bigr| = \int_t^{t(1+\eps)}
\bigl|v'(s)\bigr| \,ds \leq\bigl|v'(t)\bigr| \eps t= \psi\bigl(v(t)\bigr) t
\eps.
\]
Similarly,
\begin{eqnarray*}
\bigl|v\bigl(t(1-\eps)\bigr) - v(t)\bigr| \bck&=&\bck\int_{t(1-\eps)}^t
\bigl|v'(s\bigr)\bigr| \,ds \leq\bigl|v'\bigl(t(1-\eps)\bigr)\bigr| t\eps
\\
\bck&=&\bck\psi\bigl(v\bigl(t(1-\eps)\bigr)\bigr) t(1-\eps) \frac{\eps}{1-\eps}.
\end{eqnarray*}
Hence (\ref{cond1})
will
hold provided $\psi(v(t))t = O(v(t))$.
\end{pf}

\section{Regularity indices and consequences}
\label{Slevreg}
In this section we use the quantitative estimates obtained above
(Proposition~\ref{Tsmalltime}) to get concrete information on the
small-time behavior of $N^\La(t)$. We are particularly concerned with
power-law behavior, which as we show below turns out to be intimately
related to the notion of upper and lower indices, which arose in
seminal papers by Blumenthal and Getoor~\cite{BG} and Pruitt \cite
{pruitt} on pathwise properties of L\'evy process.

Let $X=(X_t, t\ge0)$ be a L\'evy process
with Laplace exponent $\psi$ given by (\ref{Dpsi}). We call
$\nu(dx)=x^{-2}\Lambda(dx)$ and recall that we assume that
$\Lambda(\{0\})=0$ to avoid a Kingman component. As discussed above,
we may also assume
that ${\rm supp}(\La)\subset[0,1/2)$.
The following definitions and properties of the upper-index
$\beta$ and of the lower-index $\delta$ of $X$ can be found in \cite
{BG} and \cite{pruitt}.

\begin{defn}
The \emph{upper index} is defined by
\begin{equation}
\label{Eupperind} \beta:= \inf \biggl\{ \alpha>0 \dvtx \int_{|x|\le1}
|x|^{\alpha} \nu(dx) <\infty \biggr\} \in[0,2].
\end{equation}
\end{defn}

To define the lower-index, following Pruitt \cite{pruitt}, we
introduce the function $h(x)=G(x)+K(x)+M(x)$, where
[since in our setting ${\rm supp}(\nu) \subset\R^+$ and moreover the
drift is 0]
\[
G(x)=\nu(y\dvtx y>x),\qquad K(x)=x^{-2}\int_{y\le x}
y^2\nu(dy)
\]
and
\[
M(x)=x^{-1}\biggl\llvert \int_{y\le x}
\frac{y^3}{1+y^2}\nu(dy) -\int_{y>x}\frac{y}{1+y^2}\nu(dy)
\biggr\rrvert.
\]

\begin{defn}
The \emph{lower index} is defined by
\begin{equation}
\label{Elowerind} \delta:=\inf\Bigl\{ \alpha\dvtx  \liminf_{x\to0}x^{\alpha}h(x)=0
\Bigr\}.
\end{equation}
\end{defn}

Note that the upper index $\beta$ of (\ref{Eupperind}) is similarly
given by
\[
\beta= \inf\Bigl\{ \alpha\dvtx  \limsup_{x\to0}x^{\alpha}h(x)=0
\Bigr\}.
\]
Therefore, it must be
\[
0\le\delta\le\beta\le2.
\]
The constants $\beta$ and $\delta$ characterize the asymptotic
behavior of $X$ near $0$; see (3.4) in
Pruitt \cite{pruitt} and Figure~3.
More precisely, if
$M_t:=\sup_{0 \le s \le t} |X_s|$, then
\begin{eqnarray*}
\limsup_{t \to0} M_t/t^{\kappa} &=& \cases{
 0,&\quad $\mbox{if } \kappa< 1/\beta$,
\vspace*{2pt}\cr
\infty,&\quad$\mbox{if } \kappa> 1/\beta,$}
\\
 \liminf
_{t \to0} M_t/t^{\kappa} &=& \cases{ %
0,&\quad $\mbox{if } \kappa< 1/\delta$,
\vspace*{2pt}\cr
\infty,&\quad $\mbox{if } \kappa> 1/\delta.$}
\end{eqnarray*}

In this section we show the following result:
\begin{proposition}\label{partial2}
If the lower-index $\delta$ is strictly greater than $1$, then for any
$\eps>0$,
%If the upper-index $\beta>1$, then for any $\eps>0$
\[
\frac{N^\La(t)}{t^{-1/(\beta+\eps-1)}} \to\infty \qquad\mbox{in probability,}
\]
and, for any $\eps\in(0,\delta- 1)$,
\[
\frac{N^\La(t)}{t^{-1/(\delta-\eps-1)}} \to0\qquad \mbox{in probability.}
\]
\end{proposition}

\begin{remark}
When Grey's condition for extinction holds, we know
(Lemma~\ref{psi<}) that $\beta\ge1$. However, by modifying the
construction in the \hyperref[Sexample]{Appendix},
%it is possible to find examples of L\'evy measures
%$\nu$ such that both $\beta>1$ and
%Grey's condition do not hold
it is possible to find examples such that $\beta>1$
and yet Grey's condition does not hold
(i.e., the corresponding coalescent
does not come down from infinity). See the second to last paragraph of
the \hyperref[Sexample]{Appendix}.
\end{remark}

Informally speaking, the following lemma states that as $t \to0$ the
function $q
\mapsto\psi(q)$ is of order at most $q^{\beta}$ and at least
$q^{\delta}$.
\begin{lemma} \label{psi<}
For each $\eps>0$
%jb added small enough
small enough, there exist finite constants $c_{\eps,\beta}$ and
$c_{\eps,\delta}$
such that
for all $v$ large enough
$c_{\eps,\delta} v^{\delta-\eps} \le\psi(v)\le c_{\eps,\beta}
v^{\beta+\eps}$. Hence if $\Lambda$ is such that the $\La$-coalescent
comes down from infinity, then
$\beta\geq1$.
\end{lemma}

\begin{pf}
Observe that for large $q$,
\begin{equation}
\label{controlpsi} \psi(q) \asymp q^2\int_{[0,1/q]}
x^2\nu(dx)+q\int_{[1/q,1]} x \nu(dx),\qquad  q \to\infty,
\end{equation}
where $f(q)\asymp g(q)$ means that both
$f=O(g)$ and $g=O(f)$.
Indeed, for $x \leq1/q$ one can use
Taylor's approximation to get $e^{-qx} -1 + qx \in[q^2 x^2/6,q^2
x^2/2]$ while for $x\geq1/q$ an easy computation shows $e^{-qx}
-1 + qx \in[q x/e, q x]$.

By definition (\ref{Eupperind}), we have that $\int_{[0,1]}
x^{\beta+\eps} \nu(dx) <\infty$. Therefore
\begin{eqnarray*}
%jb n-> (n+1) in the exp on the LHS
\sum_{n=0}^{\infty}
e^{-(n+1) (\beta+\eps)} \nu\bigl(\bigl[e^{-n-1},e^{-n}\bigr]\bigr) &
\le& \sum_{n=0}^{\infty} \int
_{e^{-n-1}}^{e^{-n}} x^{\beta
+\eps} \nu(dx) \\&=& \int
_{[0,1]} x^{\beta+\eps} \nu(dx)< \infty.
\end{eqnarray*}
In particular, there exists a constant $c>0$ such that for all $n\ge
1$,
\begin{equation}
\label{upboundnu} \nu\bigl(\bigl[e^{-n-1},e^{-n}\bigr]\bigr) \le c
e^{(n+1) (\beta+\eps)}.
\end{equation}
As a consequence,
%jb added the following
for $\varepsilon< 2-\beta$
\begin{eqnarray*}
\int_0^{1/q} x^2\nu(dx) &\le& \sum
_{n= \lfloor\log
q\rfloor}^{\infty} \int_{e^{-n-1}}^{e^{-n}}
x^2 \nu(dx)\\
& \le& c \sum_{n= \lfloor\log q\rfloor}^{\infty}
e^{(n+1)(\beta+\eps)}e^{-2n} \le c q^{\beta-2+\eps},
\end{eqnarray*}
where the finite positive constant $c$ may change from one inequality
to the next one.
Similarly, one estimates
\[
\int_{1/q}^1 x \nu(dx) \le c q^{\beta-1+\eps}.
\]
Together with (\ref{controlpsi}), this yields the upper bound
$
\psi(q)= O( q^{\beta+\eps})$.

For the lower bound,
recall definition (\ref{Elowerind}) and related notation.
Observe that
\[
\int_{y>x}\frac{y}{1+y^2}\nu(dy) \asymp\int
_{y>x} y\nu(dy),\qquad  x\in(0,1).
\]
The first
integral in the definition of $M(x)$ is of order $\int_{y\le x} y
\Lambda(dy) =O(x)$, so it is
negligible, in comparison.
Also, note that
%jb added
as $x \to0$,
\[
G(x)=\int_{y>x}\nu(dy) \le x^{-1}\int
_{y>x}y\nu(dy) \asymp M(x).
\]
Combining this with the definition of $K(x)$ and (\ref{controlpsi})
one gets
\begin{equation}
\label{hpsi} h(x)\asymp\psi(1/x)\qquad \mbox{as } x \to0.
\end{equation}
Due to (\ref{Elowerind}), we have $h(x)\ge Cx^{-\delta+\eps}$ for
all $x$ sufficiently small and for some $C>0$, and the lower bound for
$\psi$ now easily follows.

Finally, assume that a given $\La$-coalescent comes down from
infinity. Then by Theorem~\ref{Tmainthm}, Grey's condition \eqref
{Egrey0} is satisfied for the corresponding measure $\La$.
Since for each $\eps>0$, $\psi(q) \le cq^{\beta+ \eps}$, we deduce
that $\beta\ge1$.
\end{pf}

\begin{remark}
Note that (\ref{hpsi}) also implies the stated upper bound on~$\psi(v)$.
\end{remark}

The asymptotic behavior of $\psi(q)$ as $q \to\infty$ induces the
asymptotic behavior of $v(t)$ as $t \to0$:
\begin{coro} \label{controlv}
Assume that the $\La$-coalescent comes down from infinity.
\begin{longlist}[(ii)]
\item[(i)] If $\beta\ge1$, we have
$ \liminf_{t \to0} t^{1/(\beta-1 +\eps)} v(t)=+\infty
$ for any $\eps>0$.

\item[(ii)] If $\delta>1$, then
$ \limsup_{t \to0}
t^{1/(\delta-\eps-1)} v(t)=0$
for any $\eps\in(0,\delta-1)$.
\end{longlist}
\end{coro}

\begin{pf}
Recall (\ref{Dv}) and
let $c_{\eps,\beta}$ and $c_{\eps,\delta}$ be as in Lemma~\ref{psi<}.
It follows immediately that
%jb constants of Lemma~31 in the denominator.
\[
\frac{ v^{-(\beta+\eps) +1}}{c_{\eps,\beta}(\beta+\eps-1)} \leq \int_v^{\infty}
\frac{dq}{\psi(q)} \le \frac{ v^{-(\delta-\eps) +1}}{c_{\eps,\delta}(\delta-\eps-1)}.\vadjust{\goodbreak}
\]
Note that since $\beta\geq1$ we are able to
integrate the lower bound for $1/\psi$, but we use the
additional constraint $\delta>1$ in order
to be able to integrate the upper bound for $1/\psi$, and thus derive
the right-hand side of
the above inequality.
Setting the middle
term to $t$, after rearranging, we obtain
\[
\biggl({\frac{c_{\eps,\beta}}{\beta+\eps-1}} \biggr)^{1/(\beta
-1+\eps)} \cdot t^{-1/(\beta-1+\eps)}\le
v(t)\le \biggl(\frac{c_{\eps,\delta}}{\delta-\eps-1} \biggr)^{1/(\delta
-1-\eps)} \cdot t^{-1/(\delta-1-\eps)},
\]
implying both statements.
\end{pf}

Proposition~\ref{Tsmalltime} and Corollary~\ref{controlv} together yield Proposition~\ref{partial2} (using Theorem~1 in \cite
{bbl1} instead of Proposition~\ref{Tsmalltime} yields the same result
in the stronger almost sure sense).

\begin{appendix}
\section*{Appendix: An instructive example}
\label{Sexample}
In this appendix we discuss a class of examples that
illustrate potential difficulties in analyzing functions $\psi$ and
$v$ directly. In particular, we show that for some $\beta\neq\delta$,
one can choose the measure $\La$ in such a way that $\psi(q)$
oscillates between $q^\delta$ and $q^\beta$, resulting in analogous
oscillations for $v(t)$ between $t^{-1/(\beta-1)}$ and $t^{-1/(\delta
-1)}$. This shows that the upper and lower bounds of Proposition~\ref
{partial2} are sharp in general.
As a bonus we provide examples of $\La$-coalescents with $\delta=1$
that come down from infinity.
Let $\beta\in(1,2)$ be fixed. Set $a_n = e^{-n}$,
$n\geq0$ and for each $n\geq0$ define the interval $J_n$ as
\[
J_n = (a_{n+1},a_n].
\]
For a subsequence $(a_{n_k})_{k\geq0}$ of $(a_n)_{n\geq0}$
define the measure
\[
\nu\bigl(dx; (a_{n_k})_{k\geq0}\bigr)\equiv\nu(dx) = \sum
_k \indic {J_{n_k}}(x)
\frac{1}{x^{\beta+1}} \,dx.
\]
Then it is easy to check that for any choice of such a subsequence,
the corresponding measure $\nu$ has the upper index $\beta$.
It is moreover easily seen that if
$n_k=k$, then
[recalling (\ref{Dpsi}) and (\ref{Dv})]
$\psi(q)\asymp q^\beta$, $u(t):= \int_t^\infty \,dq/\psi(q)\asymp
(1/t)^{\beta-1}$, as
$t\to\infty$, and as a result $v(t) \asymp t^{-1/(\beta-1)}$, as
$t\to0$.
The remaining calculations however confirm that
if one chooses the intervals $J_{n_k}$ sparse
enough as $k\to\infty$, the asymptotic behavior of the functions
$\psi$, $u$ and $v$ can become quite irregular.

By (\ref{controlpsi}), estimating $u(t)$ as $t\to
\infty$ (up to constants) amounts to estimating
\[
\int_t^\infty\frac{1}{q^2 \int_{[0,1/q]} x^2 \nu(dx) + q
\int_{(1/q,1]} x\nu(dx)} \,dq.
\]
Define $k^*=k^*(q) = \max\{ k \dvtx a_{n_k} \geq1/q \}$ and set $\beta
_1:=\beta-1$ and $\beta_2:=2-\beta$
so that $\beta_1,\beta_2>0$.
First compute
\begin{eqnarray}\label{Esecondmom}
\int_{[0,1/q]} \bck x^2 \bck \nu(dx) &=&\sum
_k \int_{J_{n_k}
\cap[0,1/q]} x^{1-\beta} \,dx
\nonumber
\\
\bck&=&\bck \frac{1}{\beta_2} \sum_{k: a_{n_k} < 1/q}
\bigl(a_{n_k}^{2-\beta}- a_{n_k+1}^{2-\beta}\bigr) +
\sum_{k: a_{n_k} \geq 1/q > a_{n_k+1}} \int_{a_{n_k +1}}^{1/q}
x^{1-\beta} \,dx
\nonumber
\\[-8pt]
\\[-8pt]
\nonumber
\bck&=&\bck \frac{1}{\beta_2}\sum_{k=k^*+1}^\infty
e^{-n_k(2-\beta)} \bigl(1- e^{-(2-\beta)}\bigr) + \int_{\exp(-n_{k^*} -1) \wedge1/q}^{1/q}
x^{1-\beta} \,dx
\nonumber
\\
\bck&=&\bck \frac{1}{\beta_2} \Biggl(\bigl(1- e^{-\beta_2}\bigr)\bck\sum
_{k=k^*+1}^\infty e^{-n_k \beta_2} +
\frac{1}{q^{\beta_2}} - \frac{1}{(e^{n_{k^*}+1}\vee q)^{\beta
_2}} \Biggr),\nonumber
\end{eqnarray}
and similarly
\begin{equation}
\label{Efirstmom} \qquad\int_{(1/q,1]} x \nu(dx) = \frac{1}{\beta_1}
\Biggl( \sum_{l=1}^{k^*-1} e^{n_l \beta_1}
\bigl(e^{\beta_1}-1\bigr) + \bigl[\bigl(q \wedge e^{n_{k^*}+1}
\bigr)^{\beta_1} - e^{n_{k^*}\beta_1}\bigr] \Biggr).
\end{equation}
From now on assume $q\geq1$, and let $k^*=k^*(q)$ be as defined above.
Note that if $1/q \in J_{n_{k^*}}$
[meaning $n_{k^*} \leq \log(q) < n_{k^*}+1$], then
\begin{eqnarray*}
& &q^2 \cdot \biggl[\frac{1}{q^{\beta_2}} - \frac{1}{e^{(n_{k^*}+1)\beta_2}} \biggr]
+ q \cdot\bigl[q^{\beta_1} - e^{n_{k^*}\beta_1}\bigr]
\\
& &\qquad =q^2 \cdot \biggl[\frac{1}{q^{2-\beta}} - \frac{1}{e^{(n_{k^*}+1)(2-\beta)}} \biggr]
+ q \cdot\bigl[q^{\beta-1} - e^{n_{k^*}(\beta-1)}\bigr] \asymp
q^\beta,
\end{eqnarray*}
where for the last estimate it is best to consider separately the two cases
$\log(q) \in[n_{k^*}, n_{k^*}+1/2)$ and $\log(q) \in[n_{k^*}+1/2,
n_{k^*}+1)$.
One can check similarly that (still assuming $1/q \in J_{n_{k^*}}$)
the initial terms, corresponding to the nonnegative series
from (\ref{Esecondmom}) and (\ref{Efirstmom}), are of the order at
most $q^{\beta-2}$ and $q^{\beta-1}$, respectively. Hence,
\begin{equation}
\label{Easymptinside} \psi(q) \asymp q^\beta, \qquad\sfrac{1} {q} \in
\bigcup_k J_{n_k},
\end{equation}
which agrees well with the ``regular'' setting where $n_k=k$.
If on the contrary, $1/q \notin\bigcup_k J_{n_k}$, then $n_{l-1}+1
\leq\log(q) < n_{l}$
for $l=k^*(q)+1$, so that computations~(\ref{Esecondmom}) and
(\ref{Efirstmom}) imply
\[
\psi(q) \asymp c_1(\beta) q^2 \sum
_{k= l}^{\infty} e^{-n_k(2-\beta
)} + c_2(\beta)
q \sum_{k=1}^{l-1} e^{n_k(\beta-1)},
\]
where $c_i(\beta)\in(0,\infty)$, $i=1,2$ are constants depending
on $\beta$ only. Due to the properties of the exponential function we
then have
\begin{equation}
\label{Easymptoutside} \psi(q) \asymp q^2 e^{-n_l(2-\beta)} + q
e^{n_{l-1}(\beta-1)}, \qquad \sfrac{1} {q} \in(a_{n_l}, a_{n_{l-1}+1}].
\end{equation}
Therefore, we need to estimate up to constants, for large $t$,
\begin{eqnarray}
\label{Esumtoestim}&& \sum_k \int_{[t,\infty) \cap[e^{n_k}, e^{n_k+1})}
\frac{1}{q^\beta} \,dq
\nonumber
\\[-8pt]
\\[-8pt]
\nonumber
&&\qquad{} + \sum_l \int
_{[t,\infty) \cap
[e^{n_{l-1}+1}, e^{n_{l}})} \frac{1}{q^2 e^{-n_l(2-\beta)} + q
e^{n_{l-1}(\beta-1)} } \,dq.
\end{eqnarray}
The first series of
integrals above can easily be evaluated as being of order
\begin{equation}
\label{EsumtoestimA} \sum_{n_k \geq\log{t}} e^{-n_k (\beta-1)} \asymp
e^{-n_{k^*(t)}(\beta-1)}.
\end{equation}
Using the formula
\[
\int_a^b \frac{dx}{Bx + Cx^2} = \biggl[
\frac{1}{B}\log{\biggl\llvert \frac{Cx}{Cx+B}\biggr\rrvert }
\biggr]_a^b
\]
for each $l$ such that $\log{t} \leq
n_{l-1}+1$, the $l$th summand in the second series in (\ref
{Esumtoestim}) equals
\begin{eqnarray}
\label{EsumtoestimB} &&\frac{1}{e^{n_{l-1}(\beta-1)}}\log\biggl\llvert \frac{e^{n_l- n_{l-1}-1}
(e^{n_l(\beta-2)}e^{n_{l-1}+1}+e^{n_{l-1}(\beta-1)})}{e^{n_l(\beta
-1)}+e^{n_{l-1}(\beta-1)}} \biggr
\rrvert
\nonumber
\\[-8pt]
\\[-8pt]
\nonumber
&&\qquad\asymp\frac{(2-\beta) (n_l -
n_{l-1})}{e^{n_{l-1}(\beta-1)}},
\end{eqnarray}
since $\beta< 2$ and $n_{l-1} \leq n_l$.

Consider the following class of examples: for some $\eps\geq0$,
define inductively $m_0:=1$, $m_{r+1}:=m_r + e^{\eps m_r}$, for $r\geq
0$, and let
\[
n_{j+1}:=n_j+1\qquad \mbox{whenever } n_j
\in[m_{2r},m_{2r+1}) \mbox{ for some }r \in\N,
\]
and otherwise (here it must be $n_j=m_{2r+1}$ for some $r\in\N$) define
$n_{j+1} = m_{2r+2}=n_j+e^{\eps n_j}$.
In words, the strictly increasing sequence $(n_k)_{k\geq1}$ looks like
the simplest arithmetic progression over a long interval, then it makes
a jump (if $\eps>0$, its size is huge in comparison to the current
value of the sequence), and immediately after the sequence continues,
its slow increase by $1$ unit at a time, until the next even larger
jump, etc.

Now fix some $\eps\in(0,\beta-1)$.
Due to
(\ref{Esumtoestim})--(\ref{EsumtoestimB}), the corresponding $\La
$-coalescent comes down from infinity, since
\begin{equation}
\label{Esuminl} \sum_{l\geq1} \frac{n_l - n_{l-1}}{e^{n_{l-1}(\beta-1)}}<
\infty.
\end{equation}
Consider first the case $1/t \in\bigcup_k J_{n_k}$, and more precisely let
$1/t=e^{-n_j}=e^{-m_{2r}}$ for some $r\in\N$ (or equivalently, $\log
{t}$ is just at the beginning of the $r$th long interval where $n$
increases by increments of $1$).
Then $k^*(t)=j$ and so the expressions in~(\ref{Esumtoestim}) is of order
\begin{eqnarray*}
e^{-n_j(\beta-1)} +\sum_{l\geq j+1}\frac{(2-\beta) (n_{l} -
n_{l-1})}{e^{n_{l-1}(\beta-1)}}
&\asymp& \biggl(\frac{1}{t} \biggr)^{\beta-1} + \frac{(2-\beta)(m_{2r+2}
 - m_{2r+1})}{e^{m_{2r+1}(\beta-1)}}\\
 &\asymp&
\biggl(\frac{1}{t} \biggr)^{\beta-1}.
\end{eqnarray*}
The middle asymptotic was obtained by splitting the sum in $l$ into two
sums, one over the indices $l$ satisfying $n_{l-1}\in\bigcup_s
[m_{2s},m_{2s+1})$ and the other over the indices $l$ satisfying
$n_{l-1}\in\{m_{2s+1}\dvtx s \in\N\}$.
The first sum is easily seen to contribute another term of order
$(1/t)^{\beta-1}$, while for the second sum the dominant term is given by
the $l$ for which $n_{l-1}=m_{2r+1}$.
The final asymptotic result is obtained by noting that due to the
definition of the sequence\vspace*{1pt} $(m_s)_{s\geq1}$, we have
$(m_{2r+2} - m_{2r+1})/e^{m_{2r+1}(\beta-1)}=1/e^{m_{2r+1}(\beta
-1-\eps)}=
1/e^{(m_{2r}+e^{\eps m_{2r}})(\beta-1-\eps)}$, and
rewriting this last expression in terms of $t$ as
$
(\frac{1}{te^{t^\eps}} )^{\beta-1-\eps}=o (\frac
{1}{t} )^{\beta-1}
$,
we conclude that for $t$ of the form $t=e^{m_{2r}}$ we have
\[
u(t) = \int_t^\infty\frac{1}{\psi(q)} \,dq \asymp
\biggl(\frac
{1}{t} \biggr)^{\beta-1}
\]
as would be true for all $t$ in the regularly varying case $\eps=0$.

We now focus on the opposite case $1/t \notin\bigcup_k J_{n_k}$, and
in particular let us consider $t=e^{m_{2r+1} +1}$ for some $r\in\N$.
Suppose that $n_j$ is such that $n_j=m_{2r+1}$ and $n_{j+1} = m_{2r+1}+
e^{\eps m_{2r+1}}$. Then we have $k^*(t)=j$, and so it can be easily
checked that the contribution of (\ref{EsumtoestimA}) to $u(t)$ is
again of order $(1/t)^{\beta-1}$. However, the contribution of (\ref
{EsumtoestimB}) to $u(t)$ is of order
\begin{eqnarray*}
\sum_{l\geq j+1}\frac{(2-\beta) (n_{l} -n_{l-1})}{e^{n_{l-1}(\beta
-1)}} &\asymp&
\frac{(2-\beta)(m_{2r+2} - m_{2r+1})}{e^{m_{2r+1}(\beta-1)}}=
\frac{(2-\beta)}{e^{m_{2r+1}(\beta-1-\eps)}} \\
&\asymp&\biggl(\frac{1}{t}
\biggr)^{\beta-1-\eps}\gg \biggl(\frac{1}{t} \biggr)^{\beta-1}.
\end{eqnarray*}
So for $t$ of the form $t=e^{m_{2r+1}+1}$we have $u(t)\asymp
(1/t)^{\beta-1+\eps}$.

The above class of examples can be generalized in the following way:
instead of a fixed $\eps\in(0,\beta-1)$, one can introduce a
nonnegative sequence $(\eps_r)_{r\geq1}$, redefine
$m_0:=1$, $m_{r+1}:=m_r + e^{\eps_r m_r}$, for $r\geq0$, and keep the
old definition of $(n_j)_{j\geq1}$ in terms of $(m_r)_{r\geq1}$.

Now if $\eps_r= \beta-1$ identically for all $r$, the corresponding
coalescent does not come down from infinity, while we noted at the
beginning of the section that
the corresponding upper index is $\beta>1$.

Similarly, if $\limsup_r \eps_r = \beta-1$ (and $\eps_r<\beta-1$,
$\forall r$) where the terms close to $\beta-1$ in the sequence $(\eps
_r)_r$ are sufficiently sparse so that (\ref{Esuminl}) holds, then the
corresponding coalescent comes down from infinity. However,
(\ref{Elowerind}), (\ref{hpsi}) and (\ref{Easymptinside})--(\ref
{Easymptoutside})
imply that the lower index $\delta$ equals to $1$, while the upper
index is still $\beta>1$.
\end{appendix}
%%%%%%%%%%%%%%%%%%%%%%%%%%%%%%
%%%%%%%%%%%%%%%%%%%%%%%%%%%%%%
%%%%%%%%%%%%%%%%%%%%%%%%%%%%%%

% imsref loaded by akundreckaite, 2013-11-18 08:31:35

% zodis "Acknowledgments" paliekamas pagal autoriu

%suskaldyti doi

\printaddresses

\end{document}